\documentclass[10pt,twoside,english]{article}
\usepackage{amssymb,amsmath,babel,geometry,latexsym,graphics,tabularx,shapepar,enumerate}
\usepackage{hyperref}
\usepackage[all,2cell]{xy} \UseAllTwocells \SilentMatrices
\def\carlitz{{{\text{Car}}}}

\newtheorem{Theorem}{Theorem}

\newtheorem{Lemme}[Theorem]{Lemma}
\newtheorem{Proposition}[Theorem]{Proposition}
\newtheorem{Corollaire}[Theorem]{Corollary}

\newcommand{\ZZ}{\mathbb{Z}}
\newcommand{\FF}{\mathbb{F}}
\newcommand{\CC}{\mathbb{C}}
\newcommand{\QQ}{\mathbb{Q}}

\newcommand{\KK}{\mathbb{K}}
\newcommand{\LL}{\mathbb{L}}
\newcommand{\MM}{\mathbb{M}}
\newcommand{\NN}{\mathbb{N}}

\newcommand{\TT}{\mathbb{T}}

\newcommand{\bsb}{\boldsymbol}

\newcommand\CVD{{\hfill\hfil{\lower 2 pt\hbox{\vrule\vbox to 7pt 
{\hrule width 6pt\vfill\hrule}\vrule}}}\vskip 0.5cm}

\title{Hankel-type determinants and Drinfeld quasi-modular forms}
\author{
Vincent Bosser\footnote{LMNO, Universit\'e de Caen, Campus 2, Sciences 3,
F14032 Caen Cedex.}  \& Federico Pellarin\footnote{LaMUSE, 23, rue du Dr. Paul Michelon, 42023 Saint-Etienne Cedex.} \footnote{Both authors are supported by the contract ANR ``HAMOT", BLAN-0115-01.}}

\begin{document}

\maketitle

\hfill {{\em In the memory of David Hayes}}

\bigskip

\begin{small}
\noindent\textbf{Abstract.} In this paper we introduce a class of determinants ``of Hankel type". We use them 
to compute certain remarkable families of Drinfeld quasi-modular forms.
\end{small}

%Drinfeld quasi-modular forms} are analogues of classical quasi-modular forms
%taking values in certain complete, algebraically closed fields containing a global fields of positive characteristic.
%We introduce and study certain deformations of 
%{\em Drinfeld quasi-modular forms} by using {\em rigid analytic trivialisations} of corresponding {\em Anderson $t$-motives}.
%We show that a sub-algebra of these deformations has a natural graduation by the group $\ZZ^2\times\ZZ/(q-1)\ZZ$ and 
%an homogeneous automorphism,
%and we deduce from this and other properties {\em multiplicity estimates}. 

%The main consequence of our results on such deformations
%is a multiplicity estimate for Drinfeld quasi-modular forms. If $V=\widetilde{M}^{\leq l}_{w,m}$ denotes the 
%vector space of Drinfeld quasi-modular forms of weight $w$, type $m$ and depth $\leq l$ and if 
%$f$ is a non-zero element of it, the order of vanishing of $f$ at the ``cusp at infinity" of
%$f$ is smaller than a constant depending on $q$ multiplied by the dimension of $V$, provided that $w$ is big enough 
%depending on $q,l$.
%
%Our result seems inaccessible by dealing directly with iterative higher derivations
%on Drinfeld quasi-modular forms and requires transcendence constructions in its proof, 
%unlike classical multiplicity estimates in characteristic zero.

\medskip

\setcounter{tocdepth}{1}
\tableofcontents

\section{Introduction}
%It is well known that the $\CC$-algebra $M_\CC$ of classical modular forms for $\mathbf{SL}_2(\ZZ)$ is 
%equal to the polynomial algebra $\CC[E_4,E_6]$ 
To motivate this paper, we will first review some problem in the classical theory of 
quasi-modular forms. Let $M_\ZZ$ be the $\ZZ$-algebra 
generated by classical modular forms (for $\mathbf{SL}_2(\ZZ)$) whose $q$-expansion has coefficients in $\ZZ$. It is well known that $M_\ZZ$ is the polynomial algebra $\ZZ[E_4,E_6,\Delta]$ where $E_4,E_6$ are the normalised
(\footnote{A formal power series $\sum_{i\geq i_0}c_iq^i$ (or $\sum_{i\geq i_0}c_iu^i$) is said to be normalised 
if $c_{i_0}=1$. A modular form is normalised, by definition, if its $q$-expansion is normalised. A similar
definition will be used for quasi-modular forms and for Drinfeld quasi-modular forms.})
Eisenstein series of weights $4,6$ respectively, and where 
$\Delta=(E_4^3-E_6^2)/1728$ is the unique normalised cusp form of weight $12$, so that, in particular, $M_\ZZ$
is finitely generated.

%In other words, $M_\ZZ$ is finitely generated, isomorphic to a quotient of the polynomial ring in three variables
%$\ZZ[X,Y,Z]$ by the ideal generated by $1728Z-(X^3-Y^2)$ (in this paper, 
%we do not examine non-isobaric relations).

Let now $\widetilde{M}_\QQ$ be the $\QQ$-algebra of classical quasi-modular forms, as defined by Kaneko and Zagier in \cite{KZ}, with the additional condition that their $q$-expansions have coefficients in $\QQ$. It is easy to show that $\widetilde{M}_\QQ=\QQ[E_2,E_4,E_6]$, where $E_2$ is the (non-modular)
normalised Eisenstein series of weight $2$, so this again is a finitely generated algebra, but over $\QQ$. 

We may then formulate the following:

\medskip

\noindent{\bf Problem 1.} {\em Compute a minimal set of generators for $\widetilde{M}_\ZZ=\widetilde{M}_\QQ\cap\ZZ[[q]]$, the $\ZZ$-algebra generated
by quasi-modular forms of $\widetilde{M}_\QQ$ whose $q$-expansions have coefficients in $\ZZ$.}

\medskip

This problem is likely to be a difficult one. The examination of the $q$-expansions of the quasi-modular forms $DE_w$ with $E_w$ normalised Eisenstein series of weight $w$, $D=qd/dq$ and Clausen-von Staudt Theorem, indicate 
that the algebra $\widetilde{M}_\ZZ$ is more likely not finitely generated, in contrast with the structure of $M_\ZZ$.
What does a minimal set of generators of $\widetilde{M}_\ZZ$ look like?

In \cite{KK2}, Kaneko and Koike introduced a notion 
of {\em extremal quasi-modular form} (\footnote{Notice that in fact, the definition of extremality of Kaneko and Koike slightly differs from ours.}). An extremal quasi-modular form of weight $w$ and depth $\leq l$ is a
non-vanishing polynomial in $E_2,E_4,E_6$ which is isobaric of weight $w$, whose degree in $E_2$ is not bigger than $l$, and such that the
order of vanishing at $q=0$ of its $q$-expansion is maximal. If $w\geq 0$ is even and $l\geq 0$, such a form exists and is proportional to
a unique normalised form in $\widetilde{M}_\QQ$ denoted by $f_{l,w}$.

Kaneko and Koike, in \cite[Conjecture 2]{KK2}, made a prediction on the size of the denominators of the coefficients of such forms
which resembles in some way to a generalisation of Clausen-von Staudt Theorem. Indeed, if Conjecture 2 of loc. cit. holds, then $f_{l,w}\in\ZZ_p[[q]]$ for every prime number $p$ such that $p\geq w$, {\em provided} that $l\leq 4$.
In addition, the following question can be addressed.

\medskip

\noindent\textbf{Question. }{\em Let $l$ be a non-negative integer, and denote by $\mathcal{E}_l$ the set of $w$'s such that 
$f_{l,w}$ exists, and belongs to $\ZZ[[q]]$. For which $l$'s is $\mathcal{E}_l$ infinite?}

\medskip

Although very few of the $f_{l,w}$'s are known to have $q$-expansion defined over $\ZZ$ (\footnote{For example, it is not known whether $f_{1,14}\in\ZZ[[q]]$ but this looks true from numerical evidence.}),
the feeling that we have, after extensive numerical computations, is that $\mathcal{E}_{l}$ is infinite for $0\leq l\leq 4$ and finite for $l>4$. In these circumstances, we would
suggest to use these forms $f_{l,w}$ with $w$ in $\mathcal{E}_l$ to construct a set of generators for $\widetilde{M}_\ZZ$ but we refrain from making any kind of written prediction in this direction because this hypothesis is, so far, largely conjectural.

In this paper, we want to discuss similar problems, arising in the theory of {\em Drinfeld quasi-modular forms,}
where we have a slightly better understanding of what is going on.
Let $q=p^e$ be a power of a
prime number $p$ with $e>0$ an integer, let $\FF_q$ be the finite field with $q$ elements.
Let us consider, for an indeterminate $\theta$, the polynomial ring $A=\FF_q[\theta]$ and its fraction field $K=\FF_q(\theta)$. 

Let $K_\infty$ be the completion of $K$ for the $\theta^{-1}$-adic valuation and let us embed 
an algebraic closure of $K_\infty$ in its completion $\CC_\infty$ for the unique extension of that valuation.
Following Gekeler in \cite{Ge}, we denote by $\Omega$ the set $\CC_\infty\setminus K_\infty$, which has a structure of 
a rigid analytic space over which the group $\Gamma=\mathbf{GL}_2(A)$ acts discontinuously by homographies,
and with the usual local parameter at infinity $u$ (denoted by $t$ in \cite{Ge} and \cite{BP}).
These facts lead quite naturally to the notion of {\em Drinfeld quasi-modular forms}, rather parallel to that of classical 
quasi-modular forms for $\mathbf{SL}_2(\ZZ)$, which are studied in \cite{BP}, and to which we refer for the required background. 

Following \cite{BP}, we have 
three remarkable formal series $E,g,h\in A[[u]]$ algebraically independent over $K(u)$, representing respectively:
the $u$-expansion of a Drinfeld quasi-modular form of weight $2$, type $1$ and depth $1$ (the false Eisenstein series of weight $2$ of
Gekeler \cite{Ge}),
the $u$-expansion of an Eisenstein series of weight $q-1$ and type $0$,
and the $u$-expansion of a Poincar\'e series of weight $q+1$ and type $1$.
The first terms of these formal series are as follows, where $[i]=\theta^{q^i}-\theta$ for $i>0$ integer (see \cite[Lemma 4.2]{BP}):
\begin{eqnarray*}E&=&u+u^{q^2-2q+2}+\cdots \in u A[[u^{q-1}]]\\
g&=&1-[1]u^{q-1}-[1]u^{q^3-2q^2+2q-1}+\cdots \in A[[u^{q-1}]]\\
h&=&-u-u^{q^2-2q+2}+\cdots \in u A[[u^{q-1}]].
\end{eqnarray*}
Let $M_{w,m}$ be the $K$-vector space of {\em Drinfeld modular forms of weight $w$, type $m$,
whose $u$-expansions are defined over $K$}, which also is the space of isobaric polynomials (for weights and types)
in $g$ and $h$ with coefficients in $K$  (\footnote{Properly speaking, to call these spaces ``spaces of Drinfeld modular forms" is an abuse 
of language; these spaces are just generated by the $u$-expansions associated to such forms, but since 
we will work here with formal series in $u$ only, it looked advantageous to make the identification
between forms and formal series. We will do the same for Drinfeld quasi-modular forms; see \cite{BP} for further explanations.}).
The $K$-vector space of {\em Drinfeld quasi-modular forms
of weight $w$, type $m$ and depth $\leq l$, defined over $K$} is the space
$$\widetilde{M}^{\leq l}_{w,m}=M_{w,m}\oplus M_{w-2,m-1}E\oplus\cdots\oplus M_{w-2l,m-l}E^l.$$
All these spaces are finite dimensional subspaces of $K[[u]]$ and we may form the $K$-algebra of Drinfeld quasi-modular 
forms
$$\widetilde{M}_K=K[E,g,h]=\bigoplus_{w,m}\bigcup_l\widetilde{M}^{\leq l}_{w,m}.$$ In analogy with the Problem 1,
we have:

\medskip

\noindent{\bf Problem 2.} {\em Compute a minimal set of generators for $\widetilde{M}_A$, the $A$-algebra generated
by quasi-modular forms of $\widetilde{M}_K$ whose $u$-expansions have coefficients in $A$.}

\medskip

We say that an element $f$ of $\widetilde{M}^{\leq l}_{w,m}\setminus\{0\}$ is an {\em extremal Drinfeld quasi-modular form} if $\text{ord}_{u=0}f$ is maximal among the orders 
 at $u=0$ of non-zero elements of that vector space.
If there exists an extremal Drinfeld quasi-modular form of $\widetilde{M}^{\leq l}_{w,m}$ (\footnote{This occurs if and only if $\widetilde{M}^{\leq l}_{w,m}\neq(0)$, that is, if and only if $w\equiv 2m\pmod{q-1}$ with $w,l\geq 0$), it is unique up to multiplication by an element of $K^\times:=K\setminus\{0\}$.}), we denote by $f_{l,w,m}$ the unique normalised such form.

To present our main result, we need to define a certain double sequence of quasi-modular forms $$(E_{j,k})_{j\in\ZZ,k\geq 1}.$$
Let us denote as usual by $\Delta=-h^{q-1}\in A[[u]]$ the opposite of the unique normalised 
cusp form of weight $q^2-1$ for $\mathbf{GL}_2(A)$, and
let us extend the notation $[j]$ to non-positive integers by simply writing $[j]=\theta^{q^j}-\theta$ for $j\in\ZZ$,
so that $[0]=0$ and $[-1]=\theta^{1/q}-\theta$.
The sub-sequence $(E_{j,1})_{j\in\ZZ}$ is defined inductively in the following way. We set $E_{0,1}=E,E_{1,1}=-\frac{Eg+h}{[1]}$ and then, for $j\geq 0$, 
by 
$$E_{j+2,1}=-\frac{1}{[j+2]}(\Delta^{q^j}E_{j,1}+g^{q^{j+1}}E_{j+1,1}),$$
and for $j\leq 1$, by
$$E_{j-2,1}=-\frac{1}{\Delta^{q^{j-2}}}([j]E_{j,1}+g^{q^{j-1}}E_{j-1,1}).$$
For example, we have the following particular cases:
\begin{eqnarray*}
E_{-1,1}^q&=&-h,\\
E_{-2,1}^{q^2}&=&-hg^q,\\
E_{-3,1}^{q^3}&=&-h(g^{q+1}-[1]^qh^{q-1})^q,
\end{eqnarray*}
and in general, for all $j\leq -1$, it is possible to check that $E_{-j,1}^{q^j}$ is a Drinfeld cusp form of weight $q^j+1$ and type $1$.

Let us write
$$B_k(t):=\prod_{0\leq i<j<k}(t^{q^j}-t^{q^i})\in\FF_q[t].$$
For $k\geq 2$ and $j\in\ZZ$, we then define $E_{j,k}$ with the following {\em determinant of Hankel type}:
$$E_{j,k}=\frac{1}{B_k(\theta)}\left|\begin{array}{lllll} E_{j,1} & E_{j+1,1} & \cdots & E_{j+k-1,1}\\
E_{j-1,1}^q & E_{j,1}^q & \cdots & E_{j+k-2,1}^q\\
E_{j-2,1}^{q^2} & E_{j-1,1}^{q^2} & \cdots & E_{j+k-3,1}^{q^2}\\
\vdots & \vdots & & \vdots \\
E_{j-k+1,1}^{q^{k-1}} & E_{j-k+2,1}^{q^{k-1}} & \cdots & E_{j,1}^{q^{k-1}}
\end{array}\right|.$$

%(\footnote{A formal series $f=\sum_{i\geq i_0}c_iu^i$ with $c_{i_0}\neq0$ is said to be 
%{\em normalised} if $c_{i_0}=1$.}).

We shall show:
\begin{Theorem}\label{th2} The following properties hold, for $j\geq 0$ and $k\geq 1$.
\begin{enumerate}
\item There exists a constant $C(q,k)$ and a sequence of integers $(l_k)_{k\geq1}$ such that for all $j\geq C(q,k)$,
$$E_{j,k}\in\widetilde{M}^{\leq(q^k-1)/(q-1)}_{(q^k-1)(q^j+1)/(q-1),k}\setminus \widetilde{M}^{\leq l_{k}}_{(q^k-1)(q^j+1)/(q-1),k}$$ with $l_{k}\rightarrow\infty$ for $k\rightarrow\infty$.
\item For all $j,k$ with $j\geq 0$, we have $\text{{\em ord}}_{u=0}E_{j,k}=q^j(q^{2k}-1)/(q^2-1)$.
\item For all $j,k$ with $j\geq 0$, we have $E_{j,k}\in A[[u]]$ and $E_{j,k}$ is normalised.
\item For $k=1$ and for $k=2$ if $q\geq 3$, we have $E_{j,k}=f_{(q^k-1)/(q-1),(q^k-1)(q^j+1)/(q-1),k}$ for all $j\geq 0$.
\end{enumerate}
\end{Theorem}
In particular, for all $j,k$, $E_{j,k}$ is non-zero, property which does not seem to follow directly from the definition above.
The interest of the theorem is that it provides in an explicit way a family of {\em normalised} Drinfeld quasi-modular forms parametrised by $\ZZ_{\geq 0}\times\ZZ_{>0}$,
with {\em unbounded depths and weights}, with {\em high order of vanishing} at $u=0$, and with $u$-expansions {\em defined over $A$}. The theorem gives a partial answer to the analogue of the Question above.
Indeed, denoting by $\mathcal{E}_{l,m}$ the set whose elements are the weights $w$ such that $f_{l,w,m}$ is defined over $A$, we have the following obvious consequence of Theorem \ref{th2}.

\begin{Corollaire} If $l=1$ and for any value of $q$, or if $q\geq 3$ and $l=q+1$, we have 
that $f_{l,l(q^j+1),l}\in A[[u]]$ for all $j\geq 0$. Therefore, for the selected values of $q,l,m$, the set $\mathcal{E}_{l,m}$ has infinitely many elements.
\end{Corollaire}

It can be shown that for $k>2$, the degree of $E_{j,k}$ in $E$ is not equal to $(q^k-1)/(q-1)$, that is, it is not maximal (it is maximal only for $k=1,2$), which may 
mean that for such values, $E_{j,k}$ is not extremal. However,
the fact that $l_{k}\rightarrow\infty$ suggests that no natural threshold for the depth (as $l=4$ in the classical case,
as suggested by \cite[Conjecture 2]{KK2})
exists in the Drinfeldian framework. Moreover, the presence of infinitely many $f_{l,w,m}$'s defined over $A$ detected
by Theorem \ref{th2}
suggests that the $A$-algebra $\widetilde{M}_A$ generated by the Drinfeld quasi-modular
forms with $u$-expansions defined over $A$ could have, as a minimal set of generators, the $f_{l,w,m}$'s
with $w\in\mathcal{E}_{l,m}$ for all $l,m$'s.

\medskip

\noindent\emph{Remark.}  With the help of a formula appearing in \cite{archiv3}, it is possible to explicitly compute the $u$-expansions
of $E_{-j,1}^{q^j}$ for $j\geq0$: we have $E_{-j,1}^{q^j}=\sum_{a\in A^+}a^{q^j}u_a$ with the notations of loc. cit. These forms, which are Hecke eigenforms, are also object of investigations by A. Petrov (private communication).

\section{Determinants of Hankel's type}

An {\em inversive difference field} $(\mathcal{K},\tau)$ is the datum of a field $\mathcal{K}$ together with 
an automorphism $\tau$ that will be supposed of infinite order. The $\tau$-{\em constant subfield} $\mathcal{K}^{\tau}$ is by definition the subfield of $\mathcal{K}$ of all the elements $x\in\mathcal{K}$ such that $\tau x=x$. Every inversive difference field can be 
embedded in an {\em existentially closed} field $\mathcal{K}^{\text{ex}}$, that is a
field endowed with an extension of $\tau$ such that $\mathcal{K}^{\tau}=(\mathcal{K}^{\text{ex}})^{\tau}$,
in which every polynomial $\tau$-difference equation has at least a non-trivial solution.

We need now to choose a field $\mathcal{K}$ with {\em two} distinguished 
automorphisms to serve our purposes. Consider two indeterminates $t,u$ and the 
field of formal series $$\mathcal{R}=K((t))((u)).$$ 
The Frobenius $\FF_q$-linear endomorphism $F$ of $\mathcal{R}$ splits as a product
$$F=\chi\tau=\tau\chi,$$  
where $\chi,\tau:\mathcal{R}\rightarrow\mathcal{R}$ are respectively $K((u))$- and $\FF_q((t))$-linear,
uniquely determined by $\chi(t)=t^q$, $\tau(u)=u^q$ and $\tau \theta=\theta^q$. The perfection 
$$\mathcal{K}=\mathcal{R}^{\text{perf}}=\bigcup_{i\geq0}\FF_q(\theta^{1/q^i})((t^{1/q^i}))((u^{1/q^i}))$$
of $\mathcal{R}$ is then endowed with extensions of $\tau$ and $\chi$ such that both
the difference fields $(\mathcal{K},\tau)$ and $(\mathcal{K},\chi)$ are inversive. Also,
$\mathcal{K}^\tau$ is equal to the perfect closure $\FF_q((t))^{\text{perf}}$ of $\FF_q((t))$ in $\mathcal{K}$ and 
$\mathcal{K}^\chi$ is equal to the perfect closure $K((u))^{\text{perf}}$ of $K((u))$ in $\mathcal{K}$.

Let $x_1,\ldots,x_s$ be elements of $\mathcal{K}$. Their {\em $\tau$-wronskian} is 
the determinant:
$$W_\tau(x_1,\ldots,x_s)=\det\left(\begin{array}{cccc} x_1 & \tau x_1 & \cdots & \tau^{s-1}x_1\\
x_2 & \tau x_2 & \cdots & \tau^{s-1}x_2\\
 \vdots & \vdots & & \vdots \\
 x_s & \tau x_s & \cdots & \tau^{s-1}x_s
\end{array}\right).$$
We recall from \cite{archiv2} that $x_1,\ldots,x_s$ are $\mathcal{K}^\tau$-linearly independent if and only if $W_\tau(x_1,\ldots,x_s)\neq 0$. Similarly, the {\em $\chi$-wronskian} $W_\chi(x_1,\ldots,x_s)$ of $x_1,\ldots,x_s$ can be introduced,
and $x_1,\ldots,x_s$ are $\mathcal{K}^\chi$-linearly independent if and only if $W_\chi(x_1,\ldots,x_s)\neq 0$.

For $\bsb{f}\in\mathcal{K}$, 
we introduce the following sequence of determinants of {\em Hankel type}:
%(\footnote{A kind of $(\tau,\chi)$-difference analog of the determinant
%$\det\left(\left(\frac{\partial^{i+j}f}{\partial z_1^i\partial z_2^j}\right)_{0\leq i,j\leq k-1}\right)$
%$$\det\left(\begin{array}{ccccc}f & \frac{\partial f}{\partial z_1} & \frac{\partial^2 f}{\partial z_1^2} & \cdots & \frac{\partial^{k-1} f}{\partial z_1^{k-1}}\\ 
%\frac{\partial f}{\partial z_2}  & \frac{\partial^2 f}{\partial z_1\partial z_2}  & \frac{\partial^3 f}{\partial z_1^2\partial z_2}  & \cdots & \frac{\partial^k f}{\partial z_1^{k-1}\partial z_2}\\ 
%\vdots & \vdots & \vdots & & \vdots \\
%\frac{\partial^{k-1} f}{\partial z_2^{k-1}}  & \frac{\partial^{k} f}{\partial z_1\partial z_2^{k-1}}  & \frac{\partial^{k+1} f}{\partial z_1^2\partial z_2^{k-1}}  & \cdots & \frac{\partial^{2k-2} %f}{\partial z_1^{k-1}\partial z_2^{k-1}}
%\end{array}\right),$$
%where $f$ is a sufficiently differentiable function of two real variables.}):
$$H_k(\bsb{f})=\left|\begin{array}{ccccc}\bsb{f} & \tau\bsb{f} & \tau^2\bsb{f} & \cdots & \tau^{k-1}\bsb{f}\\ 
\chi\bsb{f} & \chi\tau\bsb{f} & \chi\tau^2\bsb{f} & \cdots & \chi\tau^{k-1}\bsb{f}\\ 
\vdots & \vdots & \vdots & & \vdots \\
\chi^{k-1}\bsb{f} & \chi^{k-1}\tau\bsb{f} & \chi^{k-1}\tau^2\bsb{f} & \cdots & \chi^{k-1}\tau^{k-1}\bsb{f}
\end{array}\right|.$$

The proposition below will be used later.

\begin{Proposition}\label{transcendence} The following conditions are equivalent.
\begin{itemize}
\item[(i)] $H_k(\bsb{f})=0$ for some $k\geq 1$.
\item[(ii)] There exist $s\geq 1$, elements $\lambda_1,\ldots,\lambda_s$ in $\FF_q((t))^{\text{perf}}$
and elements $b_1,\ldots,b_s$ in some algebraic closure $K((u))^{\text{alg}}$ of $K((u))$ such that,
in some existentially closed extension of $(\mathcal{K},\tau)$ containing $K((u))^{\text{alg}}$,
$$\bsb{f}=\lambda_1b_1+\cdots+\lambda_sb_s.$$
\item[(iii)] For some $s\geq 1$, there exist
elements $\mu_1,\ldots,\mu_s$ in $K((u))^{\text{perf}}$, and
elements $b_1',\ldots,b_s'$ in $\FF_q((t))^{\text{alg}}$,
an algebraic closure of $\FF_q((t))$, such that 
$$\bsb{f}=\mu_1b_1'+\cdots+\mu_sb_s',$$ in some existentially closed extension of $(\mathcal{K},\chi)$
containing $\FF_q((t))^{\text{alg}}$.
\end{itemize}
\end{Proposition}
\noindent\emph{Proof.} It is easy to show that each of the second and the third conditions separately implies the first.
Let us show that the first condition implies the second. 
Assuming that
$H_k(\bsb{f})=0$ for some $k\geq 1$ is equivalent to say that
$W_\chi(\bsb{f},\tau\bsb{f},\ldots,\tau^{k-1}\bsb{f})=0$. Hence, there exist $a_0,\ldots,a_s\in\mathcal{K}^{\chi}=K((u))^{\text{perf}}$ with $a_0a_s\neq 0$, such that
$$a_0\bsb{f}+a_1\tau\bsb{f}+\cdots+a_s\tau^s\bsb{f}=0.$$
On the other hand, the algebraic equation
$$a_0X+a_1X^q+\cdots+a_sX^{q^s}=0$$
has $s$ solutions $b_1,\ldots,b_s$ in an algebraic closure $K((u))^{\text{alg}}$ of $K((u))$, 
which are linearly independent over the field $(\mathcal{K}^\tau)^{F}=\FF_q$.
In particular, $W_F(b_1,\ldots,b_s)\neq0$. 

Let us consider the compositum $\mathcal{F}$ of $\mathcal{K}$ and $K((u))^{\text{alg}}$ in some
existentially closed extension of the difference field $(\mathcal{K},\tau)$ (so we embed $K((u))^{\text{alg}}$
in the existentially closed difference field $(\mathcal{K}^{\text{ex}},\tau)$).
The restriction $\tau|_{K((u))^{\text{alg}}}$ of $\tau$
is equal to the restriction of the Frobenius $F|_{K((u))^{\text{alg}}}$. Moreover, obviously,
$W_F(b_1,\ldots,b_s)=W_\tau(b_1,\ldots,b_s)$ so that $b_1,\ldots,b_s$
are also $\mathcal{F}^{\tau}$-linearly independent, $\mathcal{F}^{\tau}$ being equal to
$\FF_q((t))^{\text{perf}}$. Since $b_1,\ldots,b_s$ span the 
$\mathcal{F}^{\tau}$-vector space of solutions of the equation
$$a_0X+a_1\tau X+\cdots+a_s\tau^sX=0,$$ we obtain the second property.

The proof that the first property implies the third is similar and left to the reader, who will notice that
it suffices to transpose the matrix used to define $H_k(\bsb{f})$.\CVD

\noindent\emph{Remark.} It is easy to show, writing $H_{s,k}$ at the place of $\tau^sH_k(\bsb{f})$ for a better display, that the following formula holds:
\begin{equation}\label{Hn2}H_{s,k}^{q+1}-H_{s,k-1}^qH_{s,k+1}=H_{s-1,k}^qH_{s+1,k},\quad (s\in\ZZ,k\geq 2).
\end{equation}
Formula (\ref{Hn2}) plays a role for $(\tau,\chi)$-difference fields similar to that of Sylvester's formula expressing determinants
$\left|\left(\frac{\partial^{i+j}f}{\partial z_1^i\partial z_2^j}\right)_{0\leq i,j\leq k-1}\right|$ as in \cite{Am}.

\medskip

The elements $\bsb{f}=\sum_{i,j}c_{i,j}t^iu^j$ that we choose are either $\bsb{d}$, either $\bsb{E}=-h\tau\bsb{d}$,
where $\bsb{d}$ is the unique solution (cf. \cite{archiv}) 
in $\FF_q[t,\theta][[u]]\subset A[[t]][[u]]$
of the linear $\tau$-difference equation
\begin{equation}\label{tau-difference}
(t-\theta^q)\Delta(\tau^2X)+g(\tau X)-X=0,
\end{equation}
with $c_{0,0}=1$ and $c_{i,0}=0$ for $i>0$.
We point out that in \cite{archiv}
we have computed some coefficients of the $u$-expansion of $\bsb{d}$. See also
the remark after Lemma \ref{degree-in-t} below.

The relationship between $H_k(\bsb{d})$ and $H_k(\bsb{E})$ is simple. Since $\chi h=h$, we have
$$\chi^{i-1}\tau^{j-1}(\bsb{E})=-h^{q^{j-1}}\tau\bigl(\chi^{i-1}\tau^{j-1}(\bsb{d})\bigr)\qquad (1\leq i,j\leq k), $$
hence
\begin{equation}\label{relationship}
H_k(\bsb{E})=(-1)^k h^{1+\cdots+q^{k-1}} \tau(H_k(\bsb{d})) = (-1)^k h^{\frac{q^k-1}{q-1}} \tau(H_k(\bsb{d})).
\end{equation}
\begin{Lemme}
We have, for $j\in\ZZ$ and $k\geq 1$:
$$E_{j,k}=\left.\frac{\tau^jH_k(\bsb{E})}{B_k}\right|_{t=\theta}.$$
\end{Lemme}
\noindent\emph{Proof.} For all $k$, $H_k(\bsb{f})$ can be rewritten, thanks to the identity $\chi=F\tau^{-1}$, as
\begin{equation}\label{alternativedef}
H_k(\bsb{f})=\left|\begin{array}{ccccc}\bsb{f} & \tau\bsb{f} & \tau^2\bsb{f} & \cdots & \tau^{k-1}\bsb{f}\\ 
(\tau^{-1}\bsb{f})^q & \bsb{f}^q & (\tau\bsb{f})^q & \cdots & (\tau^{k-2}\bsb{f})^q\\ 
\vdots & \vdots & \vdots & & \vdots \\
(\tau^{1-k}\bsb{f})^{q^{k-1}} & (\tau^{2-k}\bsb{f})^{q^{k-1}} & (\tau^{3-k}\bsb{f})^{q^{k-1}} & \cdots & \bsb{f}^{q^{k-1}}\end{array}\right|.\end{equation}

It is proved in \cite{archiv} that $\bsb{E}|_{t=\theta}=E=E_{0,1}$. Moreover,
by Lemma 22 of \cite{archiv} we have  
$\tau\bsb{E}=\frac{1}{t-\theta^q}(g\bsb{E}+\bsb{h})$,
hence $(\tau\bsb{E})|_{t=\theta}=E_{1,1}$. Now, as one sees from Equation (\ref{tau-difference}) above,
or by Proposition 9 of \cite{archiv}, the function $\bsb{E}$ satisfies the 
linear $\tau$-difference equation
\begin{equation*}
(t-\theta^{q^2})(\tau^2\bsb{E})=g^q (\tau\bsb{E}) + \Delta\bsb{E}.
\end{equation*} 
It easily follows from this, by induction, that $(\tau^j\bsb{E})|_{t=\theta}$ is well defined for all $j\in\ZZ$, and
is equal to $E_{j,1}$. Comparing the definition of $E_{j,k}$ with (\ref{alternativedef}) we immediately 
recover that the forms $E_{j,k}$ of Theorem \ref{th2} are, for $k\geq 1$, precisely the formal series of $K[[u]]$
obtained by substituting
$t$ with $\theta$ in $\tau^jH_k(\bsb{E})/B_k$, a licit operation.\CVD

\section{Properties of the determinants $H_k(\bsb{d})$}

Let $k\geq 1$ be an integer. Either $H_k(\bsb{d})=0$, or there exists $\nu_k\in\ZZ_{\geq 0}$ such that
$$H_k(\bsb{d})=\sum_{s\geq\nu_k}\kappa_{k,s}u^s$$ with $\kappa_{k,s}\in\FF_q[t,\theta]$ and $\kappa_{k,\nu_k}\neq 0$. We will prove the
Theorem below, from which we will deduce Theorem \ref{th2}.
\begin{Theorem}\label{th1}
We have $H_k(\bsb{d})\neq 0$ for all $k\geq 1$, and 
the following properties hold.
\begin{enumerate}
\item $\nu_k=\frac{(q^k-1)(q^{k-1}-1)}{q^2-1}$,
\item $\kappa_{k,\nu_k}=B_k(t)$,
\item $H_k(\bsb{d})/\kappa_{k,\nu_k}$ lies in $\FF_q[t,\theta][[u^{q-1}]]$ and is normalised.
\end{enumerate}
\end{Theorem}

\begin{Corollaire}
The function $\bsb{d}$ can be expressed neither as a finite linear combination $\lambda_1b_1+\cdots+\lambda_sb_s$
with $\lambda_1,\ldots,\lambda_s\in \FF_q((t))^{\text{perf}}$ and $b_1,\ldots,b_s\in K((u))^{\text{alg}}$, nor as
a finite linear combination $\mu_1b_1'+\cdots+\mu_sb_s'$
with $\mu_1,\ldots,\mu_s\in K((u))^{\text{perf}}$ and $b_1',\ldots,b_s'\in \FF_q((t))^{\text{alg}}$.
\end{Corollaire}
\noindent\emph{Proof.} By Theorem \ref{th1}, $H_k(\bsb{d})\neq0$ for all $k$. Therefore, we can
apply Proposition \ref{transcendence}.\CVD

The rest of this section is devoted to the proof of Theorem \ref{th1}. Since the $u$-expansions of many forms involved (like $g,\Delta,\bsb{d}\ldots$),
are actually expansions in powers of $u^{q-1}$, it will be convenient to set
$$ v:= u^{q-1}.$$
In Section~\ref{subsec1}, we first prove a general divisibility property for the coefficients of the $u$-expansion of 
$H_k(f)$ for formal series $f\in\FF_q[t,\theta][[v]]$. Then, in Section~\ref{subsec2},
we carefully study the growth of the degree in $t$ of the coefficients of $\bsb{d}$. Finally,
we complete the proof of Theorem~\ref{th1} in Section~\ref{subsec3}.

%\subsection{Computing the bad reduction locus of the operators $G_k$}
\subsection{Computation of normalisation factors}\label{subsec1}

\begin{Proposition}\label{divisibility}
Let $f$ be a formal series in $\FF_q[t,\theta][[v]]$, so that we have a formal series expansion
$H_k(f)=\sum_{s\geq 0}\kappa_sv^s$ with $\kappa_s\in\FF_q[t,\theta]$ for all $s$.
Then, the  polynomial $B_k(t)$ divides $\kappa_s$ for all $s\geq 0$.
\end{Proposition}

\noindent\emph{Proof.}
We observe that if for $1\leq i,j\leq k$ we have formal expressions $f_{i,j}=\sum_{s\in\mathcal{I}}c_{i,j,s}$,
then, by multilinearity:
\begin{equation}\label{devdet}
\left|\begin{array}{ccc} f_{1,1} & \cdots & f_{1,k}\\
\vdots & & \vdots \\ f_{k,1} & \cdots & f_{k,k}\end{array}\right|=\sum_{s_1,\ldots,s_k\in\mathcal{I}}
\left|\begin{array}{ccc} c_{1,1,s_1} & \cdots & c_{1,k,s_k}\\
\vdots & & \vdots \\ c_{k,1,s_1} & \cdots & c_{k,k,s_k}\end{array}\right|.
\end{equation}

Let us write $f=\sum_{s\geq 0}c_sv^s$ with $c_s\in\FF_q[t,\theta]$. We set 
$$f_{i,j}=\chi^{i-1}\tau^{j-1}(f)=\sum_{s\geq 0}\chi^{i-1}\tau^{j-1}(c_{s}v^s)=\sum_{s\geq 0}c_{s}(t^{q^{i-1}},\theta^{q^{j-1}})u^{q^{j-1}s}$$
so that $c_{i,j,s}=c_s(t^{q^{i-1}},\theta^{q^{j-1}})v^{sq^{j-1}}$.
By (\ref{devdet}), we obtain that 
\begin{equation}\label{Hk-expansion}
H_k(f)=\sum_{s_1,s_2,\ldots,s_k}v^{s_1+s_2q+\cdots+s_kq^{k-1}}d_{s_1,s_2,\ldots,s_k},
\end{equation}
where 
\begin{equation}\label{ds-def}
d_{s_1,s_2,\ldots,s_k}=\left|\begin{array}{llll} c_{s_1}(t,\theta) & c_{s_2}(t,\theta^q) & \cdots & c_{s_k}(t,\theta^{q^{k-1}}) \\
 c_{s_1}(t^q,\theta) & c_{s_2}(t^q,\theta^{q}) & \cdots & c_{s_k}(t^q,\theta^{q^{k-1}}) \\
 \vdots & \vdots & & \vdots \\
 c_{s_1}(t^{q^{k-1}},\theta) &  c_{s_2}(t^{q^{k-1}},\theta^q) & \cdots & c_{s_k}(t^{q^{k-1}},\theta^{q^{k-1}})
\end{array}\right|.
\end{equation}

We use the fact that $c_s=\sum_{\mu}\kappa_{\mu,s}\theta^{\mu}\in\FF_q[t,\theta]$,
with $\kappa_{\mu,s}\in\FF_q[t]$. 
Let us apply (\ref{devdet}) again, this time with $f_{i,j}=\chi^{i-1}\tau^{j-1}c_{s_j}=\sum_{\mu}(\chi^{i-1}\kappa_{\mu,s_j})\theta^{\mu q^{j-1}}$
and $c_{i,j,\mu}=(\chi^{i-1}\kappa_{\mu,s_j})\theta^{\mu q^{j-1}}$.

We obtain that
$$d_{s_1,\ldots,s_k}=\sum_{\mu_1,\ldots,\mu_k}\theta^{\mu_1+\mu_2q+\cdots+\mu_{k-1}q^{k-1}}e_{\mu_1,\ldots,\mu_k},$$
where $$e_{\mu_1,\ldots,\mu_k}=\left|\begin{array}{ccc}\eta_1 & \cdots & \eta_k\\ \vdots & & \vdots \\
\chi^{k-1}\eta_1 & \cdots & \chi^{k-1}\eta_k\end{array}\right|,$$ with $\eta_j=\kappa_{\mu_j,s_j}$.
Now, by multilinearity, $e_{\mu_1,\ldots,\mu_k}$ is a sum of Moore's determinants:
$$M(\nu_1,\ldots,\nu_k)=\left|\begin{array}{ccc}t^{\nu_1} & \cdots & t^{\nu_k}\\ \vdots & & \vdots \\
t^{\nu_1 q^{k-1}} & \cdots & t^{\nu_k q^{k-1}}\end{array}\right|.$$

We then apply the following lemma, which completes the proof of Proposition \ref{divisibility}.
\begin{Lemme}\label{mitchell}
The formula
$$M(0,1,\ldots,k-1)=B_k(t)$$ holds. Moreover,
for any choice of $\nu_1,\ldots,\nu_k$, $B_k(t)$ divides $M(\nu_1,\ldots,\nu_k)$.
\end{Lemme}
\noindent\emph{Proof.} The explicit formula is a well known application, either of Moore's determinants, or Vandermonde's determinants. As for the divisibility property, this follows from an old and well known result of Mitchell, \cite{Mi},
as $M(\nu_1,\ldots,\nu_k)$ can be viewed as a generalised Vandermonde's determinant.
\CVD

\subsection{The degrees of the coefficients of $\bsb{d}$}\label{subsec2}

To prove Theorem~\ref{th1}, we will need a precise estimate of the growth of the degree in $t$ of the coefficients of $\bsb{d}$.
Recall that the function $\bsb{d}$ lies in $\FF_q[t,\theta][[v]]$, where $v=u^{q-1}$. We will write in what follows
\begin{equation}\label{eq-d}
\bsb{d}=\sum_{s\geq 0}c_s v^s,
\end{equation}
where $c_s\in A[t]$. 
The aim of this section is to prove the following lemma.
\begin{Lemme}\label{degree-in-t}
Let $s\geq 0$ and $l\geq 0$ be integers satisfying
$$s < 1+q^2+\cdots+q^{2l}.$$
Then
$$\deg_tc_s\leq l.$$
Moreover, for all $l\geq 0$ we have
\begin{equation}\label{coeff}
c_{1+q^2+\cdots+q^{2(l-1)}}(t)=(-1)^lt^l+\cdots,
\end{equation}
where the dots stand for terms of degree $<l$.
\end{Lemme}

\noindent\textbf{Remark.} We have used here the convention that the empty sum is zero, so we have $1+q^2+\cdots+q^{2(l-1)}=0$ when $l=0$. 

\noindent\emph{Proof.}
Write
\begin{equation*}
g = 1 - [1]v +\cdots = \sum_{s\geq 0}\gamma_s v^s \in A[[v]],
\end{equation*}
and
\begin{equation*}
\Delta  = -v (1-v^{q-1}+\cdots) = \sum_{s\geq 0}\delta_s v^s\in v A[[v]].
\end{equation*}
As in \cite{archiv}, we will use the following recursion formula for the coefficients $c_s$, which easily follows from
the $\tau$-difference equation (\ref{tau-difference}) (see \cite[Formula (30)]{archiv}) :
\begin{equation}\label{recursion}
c_s = \sum_{i+jq=s}\gamma_i (\tau c_j)\ + (t-\theta^q) \sum_{i+jq^2=s}\delta_i (\tau^2 c_j).
\end{equation}

We first prove by induction on $s\geq 0$ that $\deg_t c_s\leq l$  for all $l$ satisfying $1+q^2+\cdots+q^{2l}>s$.
This statement is clearly true for $s=0$ and $s=1$, since $c_0=1$ and $c_1=-(t-\theta)$.
Let now $s\geq 2$ and $l\geq 0$ be such that $s < 1+q^2+\cdots+q^{2l}$, and consider the formula (\ref{recursion}).
If $j$ is an index occurring in the first sum, then we have $j\leq s/q <s$, hence
\begin{equation}\label{eq1}
\deg_t \tau c_j= \deg_t c_j\leq l
\end{equation}
by induction hypothesis. Let now $(i,j)$ be a pair of indices occurring in the second sum.
If $i=0$, then $\delta_i=0$ and $\delta_i (\tau^2 c_j)=0$. If $i\geq 1$, then
$j\leq (s-1)/q^2 <1+\cdots+q^{2(l-1)},$
so
\begin{equation}\label{eq2}
\deg_t \tau^2c_j= \deg_t c_j\leq l-1
\end{equation}
by induction hypothesis applied to $j$ and $l-1$ (note that $j<s$). Since the coefficients $\gamma_i$ and $\delta_i$ do not depend on $t$,
it follows from (\ref{eq1}), (\ref{eq2}) and (\ref{recursion}) that $\deg_t c_s \leq l$ as required.

Let us now prove the second part of the lemma. We argue by induction on $l$. For $l=0$ and $l=1$ the assertion is true.
Let now $l\geq 2$ be an integer, and suppose that the formula (\ref{coeff}) holds for $l-1$.
Put $s:=1+\cdots+q^{2(l-1)}$, and consider again the recursion formula (\ref{recursion}).
If $j$ is any index appearing in the first sum, then, as before, $j<s=1+\cdots+q^{2(l-1)}$. Hence $\deg_t \tau c_j=\deg_t c_j\leq l-1$
by the first part of the lemma. Let us now consider a pair $(i,j)$ appearing in the second sum of (\ref{recursion}).
The smallest possible value for $i$ is $i=1$ (since $s\equiv 1 \pmod{q^2}$), for which we have
$j=1+\cdots + q^{2(l-2)}$. In this case, the induction hypothesis yields (since $\delta_1=-1$)
$$ \delta_i (\tau^2 c_j) = (-1)^{l}t^{l-1}+\cdots $$
If now $i>1$, then $j<1+\cdots + q^{2(l-2)}$, hence $\deg_t (\tau^2 c_j) = \deg_t c_j \leq l-2$ by the first part of the lemma.
It follows from these considerations that
$$ (t-\theta^q) \sum_{i+jq^2=s}\delta_i (\tau^2 c_j) = (t-\theta^q)((-1)^{l}t^{l-1}+\cdots)=(-1)^lt^l+\cdots $$
Summing up, we have proved that $c_{1+q^2+\cdots+q^{2(l-1)}}(t)=(-1)^lt^l+\cdots$ \CVD

\noindent\emph{Remark.} The following explicit formula can be deduced from (\ref{recursion}), see \cite{archiv}. 
\begin{align}\label{d-expansion}
\bsb{d} = & 1 - (t-\theta)v - (t-\theta)v^{q^2-q+1} + (t-\theta)v^{q^2} + (t-\theta)(t-2\theta^q+\theta) v^{q^2+1} \nonumber \\
& \phantom{1 - (t-\theta)v - (t-\theta)v^{q^2-q+1}}  - (t-\theta)(t-\theta^q) v^{q^2+q} + (\theta^q-\theta)(t-\theta)(t-\theta^q) v^{q^2+q+1} + \cdots
\end{align}

\subsection{Proof of Theorem~\ref{th1}}\label{subsec3}

We can now begin the proof of Theorem~\ref{th1}. We write as before
$$ \bsb{d}=\sum_{s\geq 0}c_s v^s, $$
where $c_s\in A[t]$. It will be convenient to introduce the following notation.
If $\bsb{s}=(s_1,\ldots,s_k)\in \NN^k$ (where $\NN=\{0,1,\ldots\}$), we define
$$ ||\bsb{s}||:=\sum_{i=1}^{k}s_iq^{i-1} $$
and, as in Section~\ref{subsec1}, Equation (\ref{ds-def}),
$$ d_{\bsb{s}} :=\det (\chi^{i-1}\tau^{j-1}c_{s_j})_{1\leq i,j\leq k} = \det (C_{s_1},C_{s_2},\ldots,C_{s_k}),$$
where $C_{s_j}$ is the column vector defined by
$$ C_{s_j} =
\begin{pmatrix}
\tau^{j-1}c_{s_j} \cr
\chi (\tau^{j-1}c_{s_j}) \cr
\vdots \cr
\chi^{k-1} (\tau^{j-1} c_{s_j}) \cr
\end{pmatrix}.
$$
With this notation, the formula (\ref{Hk-expansion}) writes
\begin{equation}\label{Hk}
H_k(\bsb{d})=\sum_{\bsb{s}} d_{\bsb{s}} v^{||\bsb{s}||},
\end{equation}
where ${\bsb{s}}$ runs over all $k$-tuples of $\NN^k$.
To prove Theorem~\ref{th1}, we will show that the first non zero coefficient in the $v$-expansion (\ref{Hk})
is obtained for only one multi-index $\bsb{s}$, namely for
$$ \bsb{s}_0:=(1+q^2+\cdots+q^{2(k-2)},\ldots,1+q^2,1,0).$$
This will easily yield the theorem. We will need for this three lemmas.

\begin{Lemme}\label{s0}
Set $\bsb{s}_0:=(1+q^2+\cdots+q^{2(k-2)},\ldots,1,0)\in\NN^k$. Then we have
$$ ||\bsb{s}_0|| = \frac{(q^k-1)(q^{k-1}-1)}{(q^2-1)(q-1)}$$
and
$$ d_{\bsb{s}_0}=B_k(t).$$
\end{Lemme}

\noindent\textit{Proof.}
The first part of the lemma amounts to compute the double sum
$$ \sum_{i=1}^{k}\sum_{j=0}^{k-1-i}q^{2j}q^{i-1},$$
which is an exercise left to the reader.

To prove the second part, we use Lemma~\ref{degree-in-t}, Equality (\ref{coeff}) :
$$
d_{\bsb{s}_0}=\begin{vmatrix}
(-1)^{k-1}t^{k-1}+\cdots &  & \cdots & & -t +\theta^{q^{k-2}} & 1 \\
(-1)^{k-1}t^{(k-1)q}+\cdots &  & \cdots & & -t^q +\theta^{q^{k-2}} & 1 \\
\vdots &                      & & & \vdots & \vdots \\
(-1)^{k-1}t^{(k-1)q^{k-1}}+\cdots &  & \cdots & & -t^{q^{k-1}} +\theta^{q^{k-2}} & 1
\end{vmatrix}.
$$
Let us denote by $C_1,\ldots, C_k$ the columns of this matrix. If we substract $\theta^{q^{k-2}} C_k$ to $C_{k-1}$,
then we eliminate the constant terms in $C_{k-1}$, that is, we get the new penultimate column $C_{k-1}'={}^t(-t,-t^q,\ldots,-t^{q^{k-1}})$.
By substracting now to the column $C_{k-2}$ a suitable linear combination (with coefficients in $\FF_q[\theta]$) of the last two columns, we
get the new column $C_{k-2}''={}^t(t^2,t^{2q},\ldots,t^{2q^{k-1}})$. Repeating this process for the columns $C_j$, $j=k-3,\ldots,1$,
we see by induction that
$$
d_{\bsb{s}_0}=\begin{vmatrix}
(-1)^{k-1}t^{k-1}& \cdots & -t & 1 \\
(-1)^{k-1}t^{(k-1)q} & \cdots & -t^q  & 1 \\
\vdots &                      & \vdots & \vdots \\
(-1)^{k-1}t^{(k-1)q^{k-1}} & \cdots & -t^{q^{k-1}} & 1
\end{vmatrix}.
$$
Now, this determinant is equal to
$$
\begin{vmatrix}
1 & t & \cdots & t^{k-1} \\
1 & t^q & \cdots & t^{(k-1)q} \\
\vdots &  \vdots & & \vdots \\
1 & t^{q^{k-1}} & \cdots & t^{(k-1)q^{k-1}}
\end{vmatrix},
$$
which is equal to $B_k(t)$ (Vandermonde determinant; see also Lemma~\ref{mitchell}).
\CVD

The next lemma roughly says that if a coefficient $d_{\bsb{s}}$ is not zero in (\ref{Hk}), and if we reorder the coefficients $c_{s_i}$ such that
the sequence $(\deg_t c_{s_i})_i$ is non decreasing, then the degrees $\deg_t c_{s_i}$ grow at least linearly in $i$.

\begin{Lemme}\label{permutation-degree}
Let $\bsb{d}=(s_1,\ldots,s_k)\in\NN^k$ such that $d_{\bsb{s}}\not=0$. Let $(i_1,\ldots,i_k)$ be a permutation of
the set $\{1,\ldots,k\}$ such that $$\deg_t c_{s_{i_1}}\leq \cdots \leq \deg_t c_{s_{i_k}}.$$
Then, for all $l$, we have $$\deg_t c_{s_{i_l}}\geq l-1.$$
\end{Lemme}

\noindent\textit{Proof.}
Let us write $d_{\bsb{s}}=\det(C_{s_1},\ldots,C_{s_k})$.
Suppose that there exists an $l$ such that $\deg_t c_{s_{i_l}} \leq l-2$. Then, since the operator $\tau$ does not change the degree in $t$,
the family $(c_{s_{i_1}},\ldots,\tau^{i_l-1}c_{s_{i_l}})$ consists
of $l$ polynomials in $K[t]$ of degree $\leq l-2$, so they are linearly dependent over $K$. Hence there exist
elements $\lambda_j\in K$, not all zero, such that
$$ \sum_{j=1}^{l}\lambda_j \tau^{i_j-1}c_{s_{i_j}}=0.$$
If we now apply the operator $\chi^{i-1}$ ($1\leq i\leq  k)$, we find :
$$ \sum_{j=1}^{l}\lambda_j \chi^{i-1}\tau^{i_j-1}c_{s_{i_j}}=0 \quad (1\leq i\leq  k).$$
In other words, we get $ \sum_{j=1}^{l}\lambda_j C_{i_j}=0$, that is, a non trivial linear combination of the columns $(i_1,\ldots,i_l)$ in $d_{\bsb{s}}$.
Hence $d_{\bsb{s}}=0$, which is a contradiction.
\CVD

We introduce a further notation. If $\sigma\in S_{\{1,\ldots,k\}}$ is a permutation of the set $\{1,\ldots,k\}$ and if
$\bsb{s}=(s_1,\ldots,s_k)$ is an element of $\NN^k$,
we define $\bsb{s}^{\sigma}:=(s_{\sigma(1)},\ldots,s_{\sigma(k)})$. We recall that $\bsb{s}_0$ was defined in Lemma~\ref{s0}.

\begin{Lemme}\label{permutation}
Let $\sigma$ be a permutation of the set $\{1,\ldots,k\}$ such that $\sigma\not=\mbox{\rm Id}$. Then
$$||\bsb{s}_0^{\sigma}||>||\bsb{s}_0||.$$
\end{Lemme}

\noindent\textit{Proof.}
We argue by induction on $k$. For $k=1$ there is nothing to prove. Let now $k\geq 2$ be an integer and let $\sigma$ be a permutation as in the lemma.
For $l\geq 1$, define $t_l$ by $t_l:=1+\cdots+q^{2(l-2)}$. We will also use the notation $\bsb{s}_0^{(k)}$ instead of $\bsb{s}_0$ to indicate
the dependence on $k$. Thus we have
$$\bsb{s}_0=\bsb{s}_0^{(k)}=(t_k,\ldots,t_1) \ \mbox{\rm and}\ ||\bsb{s}_0^{(k)}||=\sum_{l=1}^k t_l q^{k-l}.$$
Let further $\tau$ denote the permutation of $\{1,\ldots,k\}$ such that
$\bigl(\bsb{s}_0^{(k)}\bigr)^{\sigma}=(t_{\tau(k)},\ldots,t_{\tau(1)})$.

First, suppose that $\tau(k)=k$. Then $\tau$ induces a non trivial permutation of the set $\{1,\ldots,k-1\}$, and
$$ ||\bigl(\bsb{s}_0^{(k)}\bigr)^{\sigma}|| - ||\bsb{s}_0^{(k)}|| = \sum_{l=1}^{k-1} (t_{\tau(l)}-t_l) q^{k-l}=
q \, ||\bigl(\bsb{s}_0^{(k-1)}\bigr)^{\sigma'}|| - ||\bsb{s}_0^{(k-1)}||, $$
where $\sigma'$ is the (non trivial) permutation of $\{1,\ldots,k-1\}$ such that $\bigl(\bsb{s}_0^{(k-1)}\bigr)^{\sigma'}=(t_{\tau(k-1)},\ldots,t_{\tau(1)})$.
By induction hypothesis, it immediately follows that $||\bigl(\bsb{s}_0^{(k)}\bigr)^{\sigma}|| - ||\bsb{s}_0^{(k)}||>0$.

Suppose now that $\tau(k) \not= k$. Then
$$ ||\bigl(\bsb{s}_0^{(k)}\bigr)^{\sigma}|| \geq t_kq^{k-\tau^{-1}(k)}\geq q t_k = \frac{q(q^{2(k-1)}-1)}{q^2-1}
> \frac{(q^k-1)(q^{k-1}-1)}{(q^2-1)(q-1)},$$
hence $||\bigl(\bsb{s}_0^{(k)}\bigr)^{\sigma} > ||\bsb{s}_0^{(k)}||$ by Lemma~\ref{s0}.
\CVD

\medskip
\noindent\textit{Proof of Theorem~\ref{th1}.} We define $\bsb{s}_0=(s_{0,1},\ldots,s_{0,k})$ as in Lemma~\ref{s0}. Thus we have
$$ s_{0,l}=1+\cdots + q^{2(k-1-l)}\quad (1\leq l \leq k).$$
Let now $\bsb{s}=(s_1,\ldots,s_k)\in\NN^k$ be such that $d_{\bsb{s}}\not=0$. Choose a permutation $\sigma$ of $\{1,\ldots,k\}$
such that $\deg_t c_{s_{\sigma(k)}}\leq \cdots \leq \deg_t c_{s_{\sigma(1)}}$. By Lemma~\ref{permutation-degree} (note the different order that
we have chosen here), we have
$\deg_t c_{s_{\sigma(l)}}\geq k-l$ for all $l$. Hence, by Lemma~\ref{degree-in-t},
$$ s_{\sigma(l)}\geq 1+\cdots + q^{2(k-1-l)}=s_{0,l}, $$
or
\begin{equation}\label{sl}
s_l\geq s_{0,\sigma^{-1}(l)}.
\end{equation}
It follows, by Lemma~\ref{permutation}, that we have
$$ ||\bsb{s}|| \geq ||\bsb{s}_0^{\sigma^{-1}}|| \geq ||\bsb{s}_0||,$$
and the equality $||\bsb{s}||=||\bsb{s}_0||$ holds only if $\sigma={\rm Id}$. In that case, the inequality (\ref{sl})
shows that $||\bsb{s}||=||\bsb{s}_0||$ only if $\bsb{s}=\bsb{s}_0$. Thus, we have shown that
in the $v$-expansion (\ref{Hk}), the first non zero coefficient is $d_{\bsb{s}_0}$ :
$$ H_k(\bsb{d})=d_{\bsb{s}_0} v^{||\bsb{s}_0||}+ {\rm higher\ terms}$$
The points 1 and 2 of Theorem~\ref{th1} follow at once from this and Lemma~\ref{s0} (recall that $v=u^{q-1}$, so
$\nu_k=(q-1) ||\bsb{s}_0||$). The point 3 is then a consequence of Proposition~\ref{divisibility}. \CVD

\section{Proof of Theorem \ref{th2}}

In order to prove Theorem~\ref{th2}, we need to introduce a few notation. For any integer $l\geq 0$ and any triple
$(\mu,\nu,m)\in\ZZ\times \ZZ\times \ZZ/(q-1)\ZZ$, we denote by $\widetilde{{\cal M}}^{\leq l}_{\mu,\nu,m}$ the $K((t))$-module of almost $A$-quasimodular forms
of weight $(\mu,\nu)$, type $m$ and depth $\leq l$ (see \cite{archiv}, Section 4.2), and we set
$$\widetilde{{\cal M}}_{\mu,\nu,m}=\bigcup_{l \geq 0}\widetilde{{\cal M}}^{\leq l}_{\mu,\nu,m}.$$
We will also write $l(f)$ for the depth of the form $f$. %$(\mu(f)$, $\nu(f)$
As in \cite[\S~5.1]{archiv}, we further set $\bsb{h}=h\bsb{d}$, and we finally define
$$ \MM^{\sharp}_{\mu,\nu,m}=K((t))[g,h,\Delta^{-1}, \bsb{E}, \bsb{h}]\cap \widetilde{{\cal M}}_{\mu,\nu,m}.$$

We have :

\begin{Lemme}\label{stability}
\begin{enumerate}
\item
If $f\in \MM^{\sharp}_{\mu,\nu,m}$, then $\tau f \in \MM^{\sharp}_{q\mu,\nu,m}$ and
$\chi f\in \MM^{\sharp}_{\mu,q\nu,m}$.
\item
For all $j\geq 0$ we have $$ l(\tau^j\bsb{E})\leq 1,\quad l(\tau^j\bsb{h})\leq 1, \quad l(\chi^j\bsb{E})\leq q^j,\quad l(\chi^j\bsb{h})\leq q^j.$$
\end{enumerate}
\end{Lemme}

\noindent \textit{Proof.} We first prove that the following equalities hold :
\begin{equation}\label{formulas}
\tau\bsb{h} = \Delta \bsb{E},\quad \tau\bsb{E}=\frac{1}{t-\theta^q}(g\bsb{E}+\bsb{h}),\quad \chi\bsb{h}=(t-\theta)^q\bsb{E}^q - \frac{g}{\Delta}\bsb{h}^q,
\quad \chi\bsb{E}=\frac{\bsb{h}^q}{\Delta}.
\end{equation}
The first equality follows at once from the definitions of $\bsb{h}$ and $\bsb{E}$ and the second is Lemma~22 of \cite{archiv}.
The last one then follows from the first :
$$ \chi\bsb{E}= \chi\left(\frac{\tau\bsb{h}}{\Delta}\right)=\frac{\bsb{h}^q}{\Delta}.$$
Finally, to prove the third equality, we use the following one, which follows for instance from \cite[Proposition~9]{archiv} or
\cite[Proposition~2.7]{BP3} :
$$ \bsb{h}=\frac{t-\theta^q}{\Delta^q}(\tau^2\bsb{h}) - \frac{g}{\Delta}\tau\bsb{h}.$$
Applying $\chi$ to both sides of this equality, and using the formula $\tau\bsb{h} = \Delta \bsb{E}$, we get
$$ \chi\bsb{h} = \frac{t^q-\theta^q}{\Delta^q}\tau(\bsb{h}^q) - \frac{g}{\Delta}\bsb{h}^q = (t-\theta)^q\bsb{E}^q - \frac{g}{\Delta}\bsb{h}^q. $$
The first part of the lemma follows at once from the relations (\ref{formulas}) (we recall that $\bsb{E}\in \widetilde{{\cal M}}^{\leq 1}_{1,1,1}$ and
$\bsb{h}\in \widetilde{{\cal M}}^{\leq 0}_{q,1,1}$). The second part is a simple induction, noticing that the depth of a form $f$ in $\MM^{\sharp}_{\mu,\nu,m}$
is nothing else than the degree $\deg_{\bsb{E}}f$, when $f$ is seen as an element of the polynomial ring
$K((t))(g,h,\Delta^{-1})[\bsb{E}, \bsb{h}]$.
\CVD

We now have all the elements to prove Theorem~\ref{th2}.

\medskip
\noindent \textit{Proof of Theorem~\ref{th2}.} 
For all $i,j\in\{1,\ldots,k\}$ we have, by Lemma~\ref{stability}:
$$ \chi^{i-1}\tau^{j-1}\bsb{E}\in \widetilde{{\cal M}}^{\leq q^{i-1}}_{q^{j-1},q^{i-1},1}. $$
It follows, by a straightforward computation, that
\begin{equation}\label{equa1}
 H_k(\bsb{E})\in \widetilde{{\cal M}}^{\leq (q^k-1)/(q-1)}_{(q^k-1)/(q-1),(q^k-1)/(q-1),k}.
\end{equation}
Replacing $t$ by $\theta$, we then obtain the value of the weight and the type 
of $E_{j,k}=(\tau^jH_k(\bsb{E})/B_k)|_{t=\theta}$.
We prove the last part of the first property of the Theorem asserting that the degree in $E$ of $E_{j,k}$ is not smaller
than some integer $l_k$ with $l_k\rightarrow\infty$ as $k\rightarrow\infty$.

By the main theorem of \cite{archiv}, if $f\in\widetilde{M}^{\leq l}_{w,m}$ is non-zero and if
\begin{equation}\label{four}w\geq 4l(2q(q+2)(3+2q)l+3(q^2+1))^{3/2},\end{equation} then
$$\text{ord}_{u=0}f\leq 16q^3(3+2q)^2lw.$$ We can choose $C(q,k)$ big enough so that
if $j\geq C(q,k)$, then (\ref{four}) holds with $f=E_{j,k}$, $w=(q^k-1)(q^j+1)/(q-1)$ and $l=\deg_E(E_{j,k})$. 
Then, we get 
$$l\geq \frac{q^{j-3}}{1+q^j}\frac{1}{16(1+q)(3+2q)^2}(1+q^k)$$
so that, enlarging $C(q,k)$ if necessary, we get, for $j\geq C(q,k)$,
$$l\geq\frac{1}{32(1+q)(3+2q)^2}(1+q^k)$$ which gives the required property of growth of the sequence $(l_k)_k$.

Using now Theorem~\ref{th1} and (\ref{relationship}), we find, for all $j\geq 0$ :
\begin{equation*} %\label{equa2}
\frac{\tau^jH_k(\bsb{E})}{\kappa_{k,\nu_k}}= \frac{(-1)^k h^{q^j\frac{q^k-1}{q-1}}}{\kappa_{k,\nu_k}}\,
\tau^{j+1}\bigl(\kappa_{k,\nu_k}u^{\nu_k}+\cdots \bigr)
= (-1)^k h^{q^j\frac{q^k-1}{q-1}} (u^{q^{j+1}\nu_k}+\cdots) \in A[[t,u]].
\end{equation*}
Substituting $t=\theta$ in this equality yields
$$ E_{j,k}= (-1)^k h^{q^j\frac{q^k-1}{q-1}} (u^{q^{j+1}\nu_k}+\cdots) \in A[[u]]. $$
The properties 2, 3 of Theorem~\ref{th2} follow at once from this and from (\ref{equa1}).

It remains to show the property 4. We consider first the case $k=1$. By definition, $H_1(\bsb{E})=\bsb{E}$
and $E_{j,1}=(\tau^j\bsb{E})|_{t=\theta}$ with $\text{ord}_{u=0}E_{j,1}=q^j$. By \cite[Theorem 1.2, Proposition 2.3]{BP2}, $E_{j,1}$ is proportional to
the function $x_j$ defined there, and hence extremal. Moreover, it is normalised, so that $E_{j,1}=f_{1,q^j+1,1}$ for $j\geq 0$.

Let us assume now that $k=2,q\geq 3$. In this case, by \cite[Theorem 1.3, Proposition 2.13]{BP2},
we see that $E_{j,2}$ is proportional to the form $\xi_j$ defined there, and hence extremal. Since it is normalised and defined 
over $A$, the proof of Theorem \ref{th2} is complete.
\CVD

\end{document}